\documentclass[leqno,a4paper,12pt]{article}
\usepackage{amsmath,amsthm}

\newtheorem{Thm}{\indent Theorem}[section]
\newtheorem{Prop}[Thm]{\indent Proposition}

\newtheorem{Cor}[Thm]{\indent Corollary}

\theoremstyle{definition}
\newtheorem{Def}[Thm]{\indent Definition}
\newtheorem{Rem}[Thm]{\indent Remark}
\newtheorem{Ex}[Thm]{\indent Example}

\def\qed{{\hskip0pt\unskip\unskip\nobreak\hfil\penalty50
          \hskip1em\hbox{}\nobreak\hfil
          {\bf q.e.d.}%
          \parfillskip=0pt\finalhyphendemerits=0
          \par}\medskip}

\newenvironment{Proof}
               {{\it Proof.}\quad}
               {\qed}

\newenvironment{Proofof}[1]
               {{\it Proof of #1.}\quad}
               {\qed}

\newcommand{\Prime}{\kern3\fontdimen1\font$'$\kern-7\fontdimen1\font}

\long\def\forget#1{}

\long\def\beginSIDEREMARK#1\endSIDEREMARK
    {{\par\bigskip\advance\leftskip by 2cm
                  \advance\rightskip by -2cm\noindent
      {\bf Our own side remark:} #1
      \par\bigskip\noindent}}

\long\def\beginFORGET#1\endFORGET{#1}
\long\def\beginFORGET#1\endFORGET{}

\def\?{\ ???\ \immediate\write16{}%
\immediate\write16{Warning: There was still a question mark . . . }%
\immediate\write16{}}

\usepackage{amsmath}

\usepackage{exscale}
\usepackage{amssymb}

\newcommand{\BQ}{{\mathbb{Q}}}

\newcommand{\BZ}{{\mathbb{Z}}}

\newcommand{\CA}{{\cal A}}
\newcommand{\CB}{{\cal B}}
\newcommand{\CC}{{\cal C}}
\newcommand{\CD}{{\cal D}}

\newcommand{\Ext}{\mathop{\rm Ext}\nolimits}
\newcommand{\Hom}{\mathop{\rm Hom}\nolimits}
\newcommand{\loccit}{[loc.$\;$cit.]}

\def\halb{\frac{1}{2}}

\def\id{{\rm id}}

\newbox\mybox
\def\arrover#1{\mathrel{
       \setbox\mybox=\hbox spread 1.4em{\hfil$\scriptstyle#1$\hfil}
       \vbox{\offinterlineskip\copy\mybox
             \hbox to\wd\mybox{\rightarrowfill}}}}
\def\larrover#1{\mathrel{
       \setbox\mybox=\hbox spread 1.4em{\hfil$\scriptstyle#1$\hfil}
       \vbox{\offinterlineskip\copy\mybox
             \hbox to\wd\mybox{\leftarrowfill}}}}

\def\ontoover#1{\mathrel{
       \setbox\mybox=\hbox spread 1.4em{\hfil$\scriptstyle#1$\hfil}
       \vbox{\offinterlineskip\copy\mybox
             \hbox to\wd\mybox{\rightarrowfill\hskip-2.8mm
                               $\rightarrow$}}}}
\def\leftontoover#1{\mathrel{
       \setbox\mybox=\hbox spread 1.4em{\hfil$\scriptstyle#1$\hfil}
       \vbox{\offinterlineskip\copy\mybox
             \hbox to\wd\mybox{$\leftarrow$\hskip-2.8mm
                               \leftarrowfill}}}}
\def\longto{\longrightarrow}
\def\into{\hookrightarrow}

\def\isoto{\arrover{\sim}}

\def\longinto{\lhook\joinrel\longrightarrow}

\usepackage[curve,matrix,arrow,cmtip]{xy}
\NoComputerModernTips

\def\myxymessage{\def\messagetext
   {Here an xy-pic diagram was omitted to speed up compilation . . . }
   \immediate\write16{\messagetext}
   \hbox{\bf \messagetext}}
\def\filxymatrix#1{\myxymessage}
\def\filxyarray#1{\myxymessage}

\newdir^{ (}{{}*!/-3pt/\dir^{(}}
\newdir_{ (}{{}*!/-3pt/\dir_{(}}
\newdir^{ )}{{}*!/+3pt/\dir^{)}}
\newdir_{ )}{{}*!/+3pt/\dir_{)}}

\def\rscript#1{\hbox to 0pt{$\scriptstyle#1$\hss}}

\let\oldbullet\bullet
\def\bullet{{\mathchoice{\oldbullet}%
                        {\oldbullet}%
                        {\scriptscriptstyle\oldbullet}%
                        {\oldbullet}}}

\newcommand{\argdot}{{\;\bullet\;}}

\newcommand{\DeffgM}{\mathop{DM^{eff}_{gm}(k)}\nolimits}
\newcommand{\DeffQgM}{\mathop{\DeffgM_\BQ}\nolimits}
\newcommand{\DgM}{\mathop{DM_{gm}(k)}\nolimits}
\newcommand{\DQgM}{\mathop{\DgM_\BQ}\nolimits}
\newcommand{\Gr}{\mathop{{\rm Gr}}\nolimits}
\newcommand{\gr}{\mathop{{\rm gr}}\nolimits}
\newcommand{\TM}{\mathop{MT(k)}\nolimits}
\newcommand{\TeffM}{\mathop{MT^{eff}(k)}\nolimits}
\newcommand{\QTM}{\mathop{\TM_\BQ}\nolimits}
\newcommand{\QTeffM}{\mathop{\TeffM_\BQ}\nolimits}

\newcommand{\DQTM}{\mathop{D \! \QTM}\nolimits}
\newcommand{\DeffQTM}{\mathop{D \! \QTeffM}\nolimits}

\newcommand{\ShN}{\mathop{Shv_{Nis}(SmCor(k))}\nolimits}

\begin{document}

%

\hfuzz=3pt
\overfullrule=10pt                   


\setlength{\abovedisplayskip}{6.0pt plus 3.0pt}
\setlength{\belowdisplayskip}{6.0pt plus 3.0pt}
\setlength{\abovedisplayshortskip}{6.0pt plus 3.0pt}
\setlength{\belowdisplayshortskip}{6.0pt plus 3.0pt}

\setlength{\baselineskip}{13.0pt}
\setlength{\lineskip}{0.0pt}
\setlength{\lineskiplimit}{0.0pt}

%
%

\title{$f$-categories and Tate motives
\forget{
\footnotemark
\footnotetext{To appear in ....}
}
}
\author{\footnotesize by\\ \\
\mbox{\hskip-2cm
\begin{minipage}{6cm} \begin{center} \begin{tabular}{c}
J\"org Wildeshaus \footnote{
Partially supported by the \emph{Agence Nationale de la
Recherche}, project no.\ ANR-07-BLAN-0142 ``M\'ethodes \`a la
Voevodsky, motifs mixtes et G\'eom\'etrie d'Arakelov''. }\\[0.2cm]
\footnotesize LAGA\\[-3pt]
\footnotesize Institut Galil\'ee\\[-3pt]
\footnotesize Universit\'e Paris 13\\[-3pt]
\footnotesize Avenue Jean-Baptiste Cl\'ement\\[-3pt]
\footnotesize F-93430 Villetaneuse\\[-3pt]
\footnotesize France\\
{\footnotesize \tt wildesh@math.univ-paris13.fr}
\end{tabular} \end{center} \end{minipage}
\hskip-2cm}
\\[4cm]
}
\date{October 9, 2008}
\maketitle
\begin{abstract}
\noindent Using the theory of $f$-categories \cite[Appendix]{B},
we prove that 
the triangulated category of Tate motives over a
field $k$ is equivalent to the
bounded derived category of its heart, provided that $k$
is algebraic over $\BQ \,$. This answers 
a question asked by Levine. \\


\end{abstract}


\bigskip
\bigskip
\bigskip

\noindent {\footnotesize Math.\ Subj.\ Class.\ (2000) numbers: 
14F42 (18E10, 18E30, 19E15, 19F27). }

\eject

\tableofcontents

\bigskip


%
%

\setcounter{section}{-1}
\section{Introduction}
\label{Intro}



Let $(\CC,t)$ be a triangulated category with a $t$-structure.
Denote its heart by $\CC^0$.
It appears natural to ask the following question. \\

\emph{ ``Can the identity on $\CC^0$
be extended to a functor $D^b(\CC^0) \to \CC$ 
from the bounded derived category of $\CC^0$ to $\CC$?''} \\

To the author's knowledge, this question was first formulated
in \cite[Sect.~3.1]{BBD}. It was solved in \loccit \
under the additional hypotheses that 
(i)~$\CC$ can be embedded into the derived category $D^+(\CA)$
of complexes over an Abelian category $\CA$, which are
bounded from below, (ii)~there are enough injectives in $\CA$. 
This was then generalized in \cite[Appendix]{B},
by introducing the notion of \emph{$f$-category}.
Another generalization was developed in \cite{K}, 
where the notion of \emph{tower} over a 
category is introduced. \\

The aim of this note is to identify an easy criterion on 
the pair $(\CC,t)$ allowing for a positive answer to the above question.
According to our Theorem~\ref{1A}, it is sufficient 
that the triangulated category
$\CC$ can be embedded into the unbounded derived category over an
exact category. \\

We admit feeling somewhat uneasy about this result.
On the one hand, it is an almost immediate
consequence of any one of the approaches from \cite{B} and \cite{K}.
It is therefore hard to imagine that an expert in category theory
might find much originality in the criterion itself.
On the other hand, it does not appear to be ``well known''
in each given context where it has chances to apply.   
This latter observation is our only source of hope for 
the reader's indulgence. \\

In this note, we use $f$-categories (which are recalled in Section~\ref{2}), 
and thus follow the approach
from \cite{B}, which we think of as being the more explicit one.
The approach from \cite{K}, which yields a universal property,
and hence a more satisfactory formal aspect, would give 
the same criterion on the existence of an extension of $\id_{\CC^0}$.
We ignore whether the extensions we get
from the two approaches coincide. \\ 

We illustrate the usefulness of Theorem~\ref{1A} by an application
for which the original criterion from \cite{BBD} is insufficient.
Let $k$ be a number field, or more generally, a field which is
algebraic over $\BQ \, $. Using our criterion,  
we establish a
functor from the bounded derived category of mixed (effective)
\emph{Tate motives over $k$} \cite{L1} 
to the triangulated category of (effective) Tate motives over $k$.
Given that the induced maps on $\Ext$-groups coincide
with those which were studied in \cite{L1}, 
the results from \loccit \ then imply that
our functor is in fact an equivalence of categories
(Theorem~\ref{1E} and Corollary~\ref{1F}).
This answers 
a question asked by Levine. \\

This work was done while I was enjoying a 
\emph{modulation de service pour les porteurs de projets de recherche},
granted by the \emph{Universit{\'e} Paris~13}. 
I wish to thank D.-C.~Cisinski for 
useful discussions and comments. 


\bigskip


%
%

\section{Statement of the main results}
\label{1}



Fix a pair $(\CC,t)$ consisting of a
triangulated category $\CC$ with a $t$-structure
\cite[D\'ef.~1.3.1]{BBD}. 
The heart of the $t$-structure will be denoted by $\CC^0$.
Recall the notion of \emph{exact category} \cite[Def.\ on p.~100]{Q}.
The aim of this note is to prove the following result.

\begin{Thm} \label{1A}
Assume that $\CC$ can be embedded as a full triangulated
sub-category into $D(\CA)$, 
the derived category of (unbounded)
complexes over an exact category $\CA$. 
Fix one such embedding. \\[0.1cm]
(a)~There exists an exact functor
\[
real : D^b(\CC^0) \longto \CC 
\]
from the derived category of bounded complexes over $\CC^0$ to $\CC$,
inducing the identity on $\CC^0$, viewed as a full sub-category of
both its source and target. The functor $real$ is $t$-exact
wrt.\ the canonical $t$-structure on $D^b(\CC^0)$ and the given
$t$-structure on $\CC$. Its composition with the cohomology functor
$H : \CC \to \CC^0$ associated to $t$ equals the canonical cohomology
functor $D^b(\CC^0) \to \CC^0$. \\[0.1cm]
(b)~The functor $real$ depends functorially on the pair $(\CC,t)$
satisfying the above hypothesis on embeddability into the derived
cate\-gory of an exact cate\-gory.
More precisely, given a second triangulated category $\CD$ with
a $t$-struc\-ture, a $t$-exact functor
\[
\alpha: \CC \longto \CD
\]
inducing an exact functor $\alpha^0: \CC^0 \to \CD^0$ 
between the hearts, 
a second exact
category $\CB$, together with an exact functor $\beta: \CA \to \CB$,
and a commutative diagram
\[
\vcenter{\xymatrix@R-10pt{
        \CC \ar[d]_{\alpha} \ar@{^{ (}->}[r] &
        D(\CA) \ar[d]^{D(\beta)} \\
        \CD \ar@{^{ (}->}[r] &
        D(\CB)
\\}}
\]
of full triangulated
embeddings, then the associated functors $real$ fit into a commutative diagram
\[
\vcenter{\xymatrix@R-10pt{
        D^b(\CC^0) \ar[r]^-{real} \ar[d]_{D^b(\alpha^0)}
                  & \CC \ar[d]^{\alpha} \\
        D^b(\CD^0) \ar[r]^-{real} & \CD  
\\}}
\]
(c)~Assume that $\Hom_\CC (M,N[2]) = 0$, for any two objects $M, N$
of $\CC^0$. Then $\CC^0$ is of cohomological dimension at most one, and
\[
real : D^b(\CC^0) \longto \CC 
\]
is fully faithful. \\[0.1cm]
(d)~Assume that $\Hom_\CC (M,N[2]) = 0$, for any two objects $M, N$
of $\CC^0$. Then 
\[
real : D^b(\CC^0) \longto \CC 
\]
is an equivalence if and only if $\CC^0$ 
generates $\CC$ (as a triangulated category).
\end{Thm}

The proof of Theorem~\ref{1A} will be given in Section~\ref{2}.

\begin{Rem} \label{1B}
Note that even when $\CA$ is Abelian,
there is no hypothesis on the compatibility
of the embedding $\CC \into D(\CA)$ with the $t$-structures
on its source and target, i.e., the $t$-structure on $\CC$
is not necessarily the one induced by the canonical $t$-structure
on $D(\CA)$.   
\end{Rem}

\begin{Rem}
One of the main results of \cite{K} states that any choice of
\emph{epivalent tower} over $\CC$ determines an extension 
of the identity on $\CC^0$
\cite[Cor.~2.7]{K}.  
It seems however that in practice the existence of such a tower
is best guaranteed when $\CC$ can be embedded into the derived category
$D(\CA)$ of some exact category $\CA \, $... We have not tried to see whether
in the situation of Theorem~\ref{1A}, the approach of \cite{K}
yields the same functor $real$. 
\end{Rem}

\begin{Rem}
As follows from Theorem~\ref{1A}~(b),
the functor $real$ is not \emph{a priori} independent
of the auxiliary data given by the exact category $\CA$ and the
embedding of $\CC$ into $D(\CA)$: choose
a triangulated $t$-exact equivalence $\kappa: \CC \to \CC$ inducing
the identity on $\CC^0$. Then the functor
$real_i$ associated to a fixed embedding $i$ of $\CC$ into $D(\CA)$ and
the functor $real_{i \circ \kappa}$ associated to $i \circ \kappa$ 
satisfy the relation
\[
real_{i \circ \kappa} = \kappa^{-1} \circ real_i \; .
\]
\end{Rem}

Theorem~\ref{1A} applies in particular in the setting
of Tate motives over a field $k$ which is algebraic over 
the field $\BQ$ of rational numbers.
Recall that 
for any integer $m$, there is defined a Tate object $\BZ(m)$ in 
the category $\DgM$ of \emph{geometrical motives over $k$} 
\cite[p.~189--192]{V}. If $m \ge 0$, then $\BZ(m)$ 
belongs to the full sub-category $\DeffgM$ \cite[Thm.~4.3.1]{V}
of \emph{effective geometrical motives over $k$}. 
By imitating the construction of \loccit , replacing 
the Abelian groups of \emph{finite
correspondences} by their tensor product with $\BQ$ 
\cite[Sect.~16.2.4 and Sect.~17.1.3]{A}, one constructs the $\BQ$-linear
analogues $\DQgM$ and $\DeffQgM$ of the above categories. The images
of the Tate objects in $\DQgM$ resp.\ $\DeffQgM$ will still be denoted
by $\BZ(m)$.

\begin{Def}[cmp.\ {\cite[Def.~3.1]{L1}}]
Define the \emph{triangulated 
category of Tate motives over $k$}
as the full triangulated sub-category $\DQTM$ of $\DQgM$ generated by
the $\BZ(m)$, for $m \in \BZ$. Define the \emph{triangulated 
category of effective Tate motives over $k$}
as the full triangulated sub-category $\DeffQTM$ of $\DeffQgM$ generated by
the $\BZ(m)$, for $m \ge 0$. 
\end{Def}

According to \cite[Thm.~1.4~i), 4.2~i)]{L1}, essentially thanks to the validity
of the Beilinson--Soul\'e vanishing conjecture for number fields,
there is a canonical non-degenerate $t$-structure on $\DQTM$.
(The relation of $K$-theory of $k$
tensored with $\BQ \,$, to $\Hom_{DM_{gm}(k)_\BQ}$ is
established by work of Bloch \cite{Bl1,Bl2}; 
see \cite[Section~II.3.6.6]{L2}.)
Denote its heart by $\QTM$. This is the Abelian category of
\emph{mixed Tate motives over $k$}. It contains all Tate objects
$\BZ(m)$ \cite[Thm.~4.2~ii)]{L1}. It
is of cohomological dimension one, and
$\Hom_{\QTM} (M,N[2]) = 0$, for any two mixed Tate motives $M, N$
\cite[Cor.~4.3]{L1}. 
Applying \cite[Th.~1.4~i), ii)]{L1}
for $a = - \infty$ and $b = 0$, one sees that
the $t$-structure on $\DQTM$ induces a non-degenerate
$t$-structure on $\DeffQTM \, $, whose heart
$\QTeffM$ is the Abelian category of 
\emph{mixed effective Tate motives over $k$}. 
It contains all Tate objects $\BZ(m)$, for $m \ge 0$.
Its inclusion
into $\QTM$ induces isomorphisms of Yoneda Ext-groups.
In particular, it is also of cohomological dimension one, and
$\Hom_{\QTeffM} (M,N[2])$ vanishes, for any two mixed 
effective Tate motives $M, N$. 

\begin{Thm} \label{1E}
Let $k$ be an algebraic field extension of $\BQ$. \\[0.1cm]
(a)~There is a canonical $t$-exact functor
\[
real : D^b \bigl( \QTeffM \bigr) \longto \DeffQTM \; ,
\]
inducing the identity on $\QTeffM$. Its composition with the cohomology functor
$H : \DeffQTM \to \QTeffM$ associated to $t$ equals the canonical cohomology
functor on $D^b \bigl( \QTeffM \bigr)$. \\[0.1cm]
(b)~The functor
\[
real : D^b \bigl( \QTeffM \bigr) \longto \DeffQTM
\]
is an equivalence of triangulated categories.
\end{Thm}

\begin{Proofof}{Theorem~\ref{1E}, assuming Theorem~\ref{1A}}
Recall the definition of the category $\ShN$ 
of \emph{Nisnevich sheaves with transfers}
\cite[Def.~3.1.1]{V}.  
It is Abelian \cite[Thm.~3.1.4]{V}, and there is a canonical full
triangulated embedding 
\[
\DeffgM \longinto D^- \bigl( \ShN \bigr)
\]
into the derived category of complexes of Nisnevich sheaves
bounded from \emph{above} \cite[Thm.~3.2.6, p.~205]{V}.
Imitating the construction from \loccit \ using rational coefficients,
one shows that there is a canonical full
triangulated embedding 
\[
\DeffQgM \longinto D^- \bigl( \ShN_\BQ \bigr) \; ,
\]
where $\ShN_\BQ \subset \ShN$ denotes the full sub-cate\-go\-ry
of Nisnevich sheaves taking values in $\BQ$-vector spaces.
We thus get a canonical embedding
into $D \bigl( \ShN_\BQ \bigr)$ of any full triangulated category
$\CC$ of $\DeffQgM$. Thus, the hypothesis of Theorem~\ref{1A} is
satisfied with $\CA = \ShN_\BQ$, and a \emph{canonical} choice of
embedding for any such sub-category $\CC$,
which in addition
is equipped with a $t$-structure. This is the case in particular
for $\CC = \DeffQTM \, $. Our claim thus follows from Theorem~\ref{1A}~(a), (d):
indeed, by the results recalled before,
 $\Hom_{\DeffQTM} (M,N[2]) = 0$, for any two mixed effective Tate
motives $M, N$, and $\QTeffM$ generates $\DeffQTM$
(since $\QTeffM$ contains all $\BZ(m)$, $m \ge 0$).
\end{Proofof}

\begin{Rem}
Note that it would not be possible to perform the above proof
in the context developed in \cite[Sect.~3.1]{BBD}. 
Indeed, our triangulated category is contained in
$D^- \bigl( \ShN_\BQ \bigr)$, while \loccit \ supposes 
its immersion into some $D^+(\CA)$ \cite[p.~79]{BBD}.
\end{Rem}

\begin{Cor} \label{1F}
Let $k$ be an algebraic field extension of $\BQ$. \\[0.1cm]
(a)~There is a unique $t$-exact functor
\[
real : D^b \bigl( \QTM \bigr) \longto \DQTM 
\]
compatible with the functor from Theorem~\ref{1E},
and sending the object $\BZ(-1)$ to $\BZ(-1)$. It
induces the identity on $\QTM$. Its composition with the cohomology functor
$H : \DQTM \to \QTM$ associated to $t$ equals the canonical cohomology
functor on $D^b \bigl( \QTM \bigr)$. \\[0.1cm]
(b)~The functor
\[
real : D^b \bigl( \QTM \bigr) \longto \DQTM
\]
is an equivalence of triangulated categories.
\end{Cor}

\begin{Proof}
By \cite[Thm.~4.3.1]{V}, the object $\BZ(1)$ is quasi-invertible 
in the cate\-gory $\DeffQTM \, $. Thus, the category $\DQTM$ is canonically 
equi\-valent
to the category obtained from $\DeffQTM$ by inverting $\BZ(1)$.

Since $\BZ(1)$ is invertible in $D^b \bigl( \QTM \bigr)$, the functor
from Theorem~\ref{1E}~(a) extends uniquely to a functor
\[
real : D^b \bigl( \QTM \bigr) \longto \DQTM 
\]
mapping $\BZ(-1)$ to $\BZ(-1)$. It is easily seen to have all the
properties listed in part~(a) of our claim.

In order to prove part~(b), note that by \cite[Thm.~4.3.1]{V} again, and
by Theorem~\ref{1E}~(b), the object $\BZ(1)$ is quasi-invertible 
in $D^b \bigl( \QTeffM \bigr)$. Therefore, the canonical
functor from the category obtained from $D^b \bigl( \QTeffM \bigr)$
by inverting $\BZ(1)$ to $D^b \bigl( \QTM \bigr)$
is an equivalence. Hence so is $real$.
\end{Proof}

\begin{Rem}
Corollary~\ref{1F} gives an affirmative answer to the question
asked in \cite[top of p.~237]{L2}, and provides a proof of 
\cite[Prop.~20.2.3.1]{A}. 
\end{Rem}

\begin{Rem}
As the proofs show,
parts~(a) of Theorem~\ref{1E} and Corollary~\ref{1F} 
hold more generally for any perfect base field
$k$ satisfying the Beilin\-son--Soul\'e vanishing conjecture (i.e.,
such that the data from \cite{L1} define a $t$-structure on $\DQTM$).
\end{Rem}


\bigskip


%
%

\section{Review of $f$-categories}
\label{2}



In this section, we review definitions and
results of Beilinson's paper \cite{B}.
An additional, rather elementary observation (Proposition~\ref{2Main})
will then allow to give a proof of Theorem~\ref{1A}. 

\begin{Def}[{\cite[Def.~A~1]{B}}] \label{2A}
(a)~An \emph{$f$-category} (or \emph{filtered triangulated category})
is a triangulated category $\CC F$, 
together with a quadruple $f = (\CC F~(\le 0) , \CC F~(\ge 0) , s , \iota)$,
whose first and second components are full 
triangulated sub-categories of $\CC F$
closed under isomorphisms in $\CC F$, $s$ is an exact auto-equivalence
on $\CC F$, and $\iota$
a transformation of functors from the identity on $\CC F$ to $s$, 
such that, putting
\[
\CC F~(\le n) := s^n \CC F~(\le 0) \quad , \quad
\CC F~(\ge n) := s^n \CC F~(\ge 0) \quad \forall \; n \in \BZ \; ,
\]
the following conditions are satisfied.
\begin{enumerate}
\item[(1)] 
We have the inclusions
\[
\CC F~(\le 0) \subset \CC F~(\le 1) \quad , \quad
\CC F~(\ge 0) \supset \CC F~(\ge 1) \; ,
\]
and the equality
\[
\bigcup_{n \in \BZ} \CC F~(\le n) = \CC F = \bigcup_{n \in \BZ} \CC F~(\ge n)
\]
(on objects).
\item[(2)]
For any object $X$ in $\CC F$, we have
\[
\iota_X = s(\iota_{s^{-1} X}) : X \longto sX \; .
\]  
\item[(3)] 
For any pair of objects $X \in \CC F~(\ge 1)$ and $Y \in \CC F~(\le 0)$,
we have
\[
\Hom_{\CC F}(X,Y) = 0 \; ,
\]
and the maps 
\[
\iota_* : \Hom_{\CC F}(Y,s^{-1}X) \longto \Hom_{\CC F}(Y,X)
\]
and
\[
\iota^* : \Hom_{\CC F}(sY,X) \longto \Hom_{\CC F}(Y,X)
\]
are both isomorphisms.
\item[(4)] 
For any object $X \in \CC F$, there exists an exact triangle
\[
A \longto X \longto B \longto A[1]
\]
in $\CC F$, such that $A \in \CC F~(\ge 1)$ and $B \in \CC F~(\le 0)$. 
\end{enumerate}
(b)~An \emph{$f$-functor} between $f$-categories is an exact functor
that respects the quadruple $f$, i.e., that conserves the sub-categories
$\CC F~(\le 0)$ and $\CC F~(\ge 0)$, and commutes with $s$ 
and $\iota$. \\[0.1cm]
(c)~Let $\CC$ be a triangulated category. An \emph{$f$-category over $\CC$}
is a pair $(\CC F , \jmath)$ consisting of an $f$-category $\CC F$ and
an equivalence of triangulated categories 
\[
\jmath : \CC \isoto \CC F~(\le 0) \cap \CC F~(\ge 0) \; .
\] 
Here, the target $\CC F~(\le 0) \cap \CC F~(\ge 0)$ denotes the
full triangulated sub-category of $\CC F$ whose objects lie both
in $\CC F~(\le 0)$ and in $\CC F~(\ge 0)$.
\end{Def}

\begin{Ex}[{\cite[Ex.~A~2]{B}}] \label{2B}
Fix an Abelian category $\CA$, and consider 
the \emph{filtered derived category} $D^\bullet F(\CA)$,
for $\bullet \in \{ \emptyset , b , - , + \}$,
over $\CA$ \cite[Sect.~3.1.1]{BBD}.
Recall that by definition, this category is obtained from the category of
complexes $A^*$
over $\CA$ equipped with a \emph{finite} decreasing filtration $F^\bullet$
by sub-complexes of $A^*$.
One sets
\[
D^\bullet F(\CA)~(\le 0) :=  \{ (A^*,F^\bullet A^*) \, , \, 
\Gr_{F^\bullet}^i A^* = 0 \; \forall \; i > 0 \}
\]
and
\[
D^\bullet F(\CA)~(\ge 0) :=  \{ (A^*,F^\bullet A^*) \, , \, 
\Gr_{F^\bullet}^i A^* = 0 \; \forall \; i < 0 \} \; . 
\]
Thus, the category $D^\bullet F(\CA)~(\le 0)$ is obtained from
the category of complexes over $\CA$ equipped with a finite decreasing 
filtration concentrated in non-positive degrees, and similarly for
$D^\bullet F(\CA)~(\ge 0)$.
The auto-equivalence $s$ is given by the shift of the filtration:
\[
s: (A^*,F^\bullet A^*) \longmapsto (A^*,F^{\bullet - 1} A^*) \; .
\]
Thus for example, an object of $D^\bullet F(\CA)$ lies
in the image $D^\bullet F(\CA)~(\le 1)$ of 
$D^\bullet F(\CA)~(\le 0)$ under $s$ if and only if 
its filtration is concentrated in degrees $\le 1$.
The natural transformation $\iota: \id \to s$ is given by the identity
on the underlying complexes. Altogether, the data
\[
(D^\bullet F(\CA)~(\le 0) , D^\bullet F(\CA)~(\ge 0) , s , \iota)
\]
define a structure of $f$-category on $D^\bullet F(\CA)$.
In fact, $D^\bullet F(\CA)$ is an $f$-category over $D^\bullet (\CA)$:
the equivalence
\[
\jmath : D^\bullet (\CA) \isoto 
D^\bullet F(\CA)~(\le 0) \cap D^\bullet F(\CA)~(\ge 0)
\]
is induced by $A^* \mapsto (A^*, Tr^\bullet)$,
where $Tr^\bullet$ is the trivial filtration concentrated in degree zero.
\end{Ex}

\begin{Rem} \label{2Bb}
We leave it to the reader to show
that the construction described in Example~\ref{2B}
generalizes to give an $f$-category $D^\bullet F(\CA)$ over $D^\bullet (\CA)$,
for any exact category $\CA$.
\end{Rem}

Beilinson proved that the forgetful functor from $D^\bullet F(\CA)$
to $D^\bullet (\CA)$ generalizes to the context of $f$-categories.

\begin{Prop}[{\cite[Prop.~A~3~(iii)]{B}}] \label{2Ba}
Let $(\CC F,\jmath)$ be an $f$-category over $\CC$.
Then there is a unique exact functor
\[
\omega: \CC F \longto \CC
\]
satisfying the following properties (1)--(3).
\begin{enumerate}
\item[(1)]~The restriction of $\omega$ to $\CC F~(\le 0)$ is left adjoint
to the embedding
\[
\CC \stackrel{\jmath}{\longinto} \CC F~(\le 0) \cap \CC F~(\ge 0)
\longinto \CC F~(\le 0) \; .
\]
\item[(2)]~The restriction of $\omega$ to $\CC F~(\ge 0)$ is right adjoint
to the embedding
\[
\CC \stackrel{\jmath}{\longinto} \CC F~(\le 0) \cap \CC F~(\ge 0)
\longinto \CC F~(\ge 0) \; .
\]
\item[(3)]~For any object $X$ of $\CC F$, the morphism
\[
\omega(\iota_X) : \omega(X) \longto \omega(sX)
\]
is an isomorphism.
\end{enumerate}
The functor $\omega$ also satisfies the following property.
\begin{enumerate}
\item[(4)]~For any pair of objects $X \in \CC F~(\le 0)$ and 
$Y \in \CC F~(\ge 0)$,
the map 
\[
\omega : \Hom_{\CC F}(X,Y) \longto \Hom_{\CC}(\omega X , \omega Y)
\]
is an isomorphism.
\end{enumerate}
\end{Prop}

\begin{Def}[{\cite[Def.~A~4]{B}}] \label{2C}
Let $(\CC F,\jmath)$ be an $f$-category over $\CC$.
Assume given $t$-structures $(\CC^{t \le 0} , \CC^{t \ge 0})$
and $(\CC F^{t \le 0} , \CC F^{t \ge 0})$ on $\CC$ and on $\CC F$.
They are said to be \emph{compatible} with each other
if $\jmath: \CC \into \CC F$
is $t$-exact (in other words, the $t$-structure on $\CC$
is induced by the $t$-structure on $\CC F$ via $\jmath$), and
\[
s (\CC F^{t \le 0}) = \CC F^{t \le -1} \; .
\]
\end{Def}

Here is the reason why we are interested in $f$-categories.

\begin{Thm}[{\cite[Prop.~A~5, Sect.~A~6]{B}}] \label{2D} 
Let $(\CC F,\jmath)$ be an $f$-cate\-gory over $\CC$.
Assume given a $t$-structure on $\CC$. 
Denote its heart by $\CC^0$. \\[0.1cm]
(a)~There is a unique $t$-structure on $\CC F$ compatible with
the $t$-structure on $\CC$. Denote its heart by $\CC F^0$. \\[0.1cm]
(b)~There is a canonical equivalence of categories
\[
\eta: \CC F^0 \isoto C^b(\CC^0)
\]
between $\CC F^0$ 
and the category of bounded complexes over $\CC^0$. \\[0.1cm]
(c)~The composition
\[
\widetilde{real}: C^b(\CC^0) \stackrel{\eta^{-1}}{\longto} \CC F^0
\longinto \CC F \stackrel{\omega}{\longto} \CC
\]
factors uniquely through an exact functor
\[
real : D^b(\CC^0) \longto \CC \; .
\]
The functor $real$ induces the identity on $\CC^0$, and is $t$-exact.
Its composition with the cohomology functor
$H : \CC \to \CC^0$ associated to $t$ equals the canonical cohomology
functor $D^b(\CC^0) \to \CC^0$.
\end{Thm}

Our contribution to the abstract theory
of $f$-structures reads as follows. 

\begin{Prop} \label{2Main}
Let $(\CC F,\jmath)$ be an $f$-cate\-gory over $\CC$,
and $\CD$ a full triangulated sub-category of $\CC$.
Then there is a unique sub-$f$-category $\CD F$ of $\CC F$
such that $\CD F$, together with
the restriction of $\jmath$ to $\CD$, forms an $f$-category 
over $\CD$.
\end{Prop}

\begin{Proof}
We need to recall a last construction due to Beilinson.
According to \cite[Prop.~A~3~(i)]{B}, there are exact functors
\[
\sigma_{\le n} : \CC F \longto \CC F~(\le n)
\]
and  
\[
\sigma_{\ge n} : \CC F \longto \CC F~(\ge n)
\]
which are left resp.\ right adjoint to the inclusions, for all integers $n$.
The adjunction properties formally imply the relations
\[
s \circ \sigma_{\le n} = \sigma_{\le n+1} \circ s
\] 
and 
\[
s \circ \sigma_{\ge n} = \sigma_{\ge n+1} \circ s
\] 
for all $n$. Again by \cite[Prop.~A~3~(i)]{B},
the functors $\sigma_{\le n}$ and $\sigma_{\ge n}$ 
respect all sub-categories $\CC F~(\le m)$ and $\CC F~(\ge m)$,
and there are canonical isomorphisms
\[
\sigma_{\le n} \circ \sigma_{\ge m} \isoto 
\sigma_{\ge m} \circ \sigma_{\le n} \; .
\]
The case $m = n$ will be of particular interest.
We identify $\sigma_{\le n} \circ \sigma_{\ge n}$
and $\sigma_{\ge n} \circ \sigma_{\le n}$, and 
note that its target equals the category $\CC F~(\le n) \cap \CC F~(\ge n)$. 
Define the exact functor 
\[
\gr_f^n : \CC F \longto \CC
\]
as the composition of $\sigma_{\le n} \circ \sigma_{\ge n}$,
of $s^{-n}$, and of $\jmath^{-1}$. Note the relation
\[
\gr_f^n = \jmath^{-1} \circ \sigma_{\le 0} 
\circ \sigma_{\ge 0} \circ s^{-n} \; .
\]
By \cite[Prop.~A~3~(ii)]{B}, for any integer $n$, 
and any object $X$ of $\CC F$, the adjunction morphisms
fit into an exact triangle
\[
\sigma_{\ge n+1} X \longto X \longto \sigma_{\le n} X \longto 
\sigma_{\ge n+1} X[1] \; ,
\]
which is unique up to unique isomorphism.

In this precise sense, any object of $\CC F$ is a successive
extension of objects of $\CC$. The analogous statement must be true for
the $f$-category $\CD F$ over $\CD$ we intend to construct.
Conversely, any successive extension of objects of $\CD$ in $\CC F$ must
belong to $\CD F$.

Hence our only choice is to
define $\CD F$ as the full sub-category of $\CC F$
of objects $X$ such that $\gr_f^n X$ belongs to $\CD \subset \CC$,
for all $n$. All $\gr_f^n$ being exact, and $\CD$ triangulated,
the category $\CD F$ is triangulated. 
The quadruple $(\CD F~(\le 0) , \CD F~(\ge 0) , s , \iota)$
is induced from $(\CC F~(\le 0) , \CC F~(\ge 0) , s , \iota)$.
More precisely, we define
\[
\CD F~(\le 0) := \CD F \cap \CC F~(\le 0)
\]
and 
\[
\CD F~(\ge 0) := \CD F \cap \CC F~(\ge 0) \; ,
\]
and $s$ and $\iota$ as the restrictions from 
the corresponding data on $\CC F$. Indeed, $s$ and $s^{-1}$ respect
$\CD F$ since $\gr_f^n \circ s = \gr_f^{n-1}$ for all integers $n$.

First, in order to verify conditions~\ref{2A}~(a)~(1)--(4), observe that 
(1)--(3) are obvious by construction. Condition~(4) is equivalent to
stating that the functors $\sigma_{\le 0}$ and $\sigma_{\ge 1}$
respect $\CD F$. In order to see this, it suffices to note that the
composition $\gr_f^n \circ \sigma_{\le 0}$ equals $\gr_f^n$ if $n \le 0$,
and zero otherwise, and similarly for $\gr_f^n \circ \sigma_{\ge 1}$.

Now consider the triangulated category $\CD F~(\le 0) \cap \CD F~(\ge 0)$.
Via $\jmath^{-1}$, it is equivalent to a full triangulated sub-category
of $\CC$. Since all $\gr_f^n$ are trivial on $\CD F~(\le 0) \cap \CD F~(\ge 0)$
except for $n = 0$, this sub-category equals $\CD$.
Therefore, $\CD F$ is indeed an $f$-category over $\CD$.
\end{Proof}

\begin{Cor} \label{2E}
Let $\CA$ be an exact category, and 
$\CC$ a full triangulated sub-category of $D(\CA)$. 
Then the embedding of $\CC$ into $D(\CA)$ induces a 
choice of $f$-category over $\CC$.
\end{Cor}

\begin{Proof}
Use Proposition~\ref{2Main} and Remark~\ref{2Bb}.
\end{Proof}

\medskip

\begin{Proofof}{Theorem~\ref{1A}}
We are in the situation of Corollary~\ref{2E},
and thus get an $f$-category $\CC F$ over $\CC$.
In addition, a $t$-structure on $\CC$ is given.
We need to equip $\CC F$ with a $t$-structure, and take the one from
Theorem~\ref{2D}~(a). Then Theorem~\ref{2D}~(c) gives the
construction of the functor $real$ and states that it has the
properties listed in part~(a) of our claim.

Part~(b) is a special case of the functorial behaviour of $real$
under $f$-functors \cite[Lemma~A~7.1]{B}. 

To prove parts~(c) and (d), note that our functor
\[
real : D^b(\CC^0) \longto \CC 
\]
is fully faithful if and only if for any two objects
$M, N$ of $\CC^0$, and any integer $p \ge 0$, the morphism
\[
\Ext^p_{\CC^0} (M,N) \longto \Hom_\CC (M,N[p])
\]
($\Ext^p = $ Yoneda Ext-group of $p$-extensions)
induced by $real$ is an isomorphism. This is obviously true for $p = 0$.
Now recall that the Yoneda Ext-groups form a universal
$\delta$-functor \cite[Prop.~4.1, Prop.~4.3]{Bu}. 
Therefore, the above morphism 
equals the value on $M$ of the $p$-th member of the natu\-ral
transformation of functors 
\[
\tau^p: \Ext^p_{\CC^0} (\argdot,N) \longto \Hom_\CC (\argdot,N[p])
\]
corresponding to this universal property.
Using the criterion from \cite[Prop.~4.2]{Bu}, one shows abstractly
(see for example \cite[p.~3]{DG})
that $\tau^p$ is an isomorphism
for $p=1$, and injective for $p=2$.  
Our hypothesis on the vanishing of $\Hom_\CC (\bullet,N[2])$
thus trivially implies that $\tau^2$ is an isomorphism.
It also implies that both source and target of $\tau^p$
are zero for all $p \ge 2$. 
\end{Proofof}


\bigskip

%
%

\end{document}